%% file: RSVarxiv.tex
\documentclass[11pt]{amsart} 
\usepackage{amsfonts}
\usepackage{amsmath}
\usepackage{amssymb}
\usepackage{amsthm}
\usepackage{cite}
\usepackage{epsfig} 
\usepackage{rotate}
\usepackage{tikz}
\usepackage{pgfplots}
\usetikzlibrary{math}

\newtheorem{theorem}{Theorem} 
\newtheorem{corollary}[theorem]{Corollary}
\newtheorem{lemma}[theorem]{Lemma}

\def\Q{\mathbb Q}
\def\Z{\mathbb{Z}}

\def\a{a}
\def\b{b}

\def\bi{{\mathrm{\bf i}}}
\def\bj{{\mathrm{\bf j}}}
\def\rg{\mathrm{g}}
\def\ri{{\mathrm{i}}}
\def\rj{{\mathrm{j}}}
\def\rij{(\ri,\rj)}
\def\rz{\mathrm{z}}
\def\cC{\mathcal{C}}
\def\cE{\mathcal{B}}
\def\cF{\mathcal{F}}
\def\cL{\mathcal{L}}
\def\cN{\mathcal{N}}
\def\cLij{\cL_{\ri,\rj}}
\def\rF{\mathrm{F}}
\def\rg{\mathrm{g}}
\def\z{\zeta}
\def\zu{\zeta\hskip-2.5pt\uparrow}
\def\zd{\zeta\hskip-2.5pt\downarrow}
\def\zud{\zeta\hskip-3.5pt\uparrow\hskip-1.5pt\downarrow}
\def\gI{I}
\def\cB{\mathcal{B}}

\def\n2{\mathbf{k}}
\def\I{\mathrm{I}}
\def\II{\mathrm{II}}
\def\III{\mathrm{III}}
\def\IV{\mathrm{IV}}

\def\qad{\hskip 5pt}
\def\mod#1{\,({\rm mod\ }#1) }
\def\proof{\noindent {\sc Proof.}\hskip 10pt}
\def\endproof{\hfill\vbox{\hrule \hbox{\vrule\kern4pt\vbox{\kern4pt
\kern4pt}\kern4pt\vrule}\hrule}\smallskip}
\def\flo#1{\lfloor #1\rfloor}
\def\cei#1{\lceil #1\rceil}

\begin{document}
\title[]{Critical curves of rotations}
\author{JAG~Roberts} 
\address{School of Mathematics and Statistics,
University of New South Wales,
Sydney, NSW 2052, Australia}
\email{jag.roberts@unsw.edu.au}
%
\author{Asaki Saito}
\address{Future University Hakodate,
116--2 Kamedanakano-cho,
Hakodate, Hokkaido 041--8655, Japan}
\email{saito@fun.ac.jp}
\author{Franco Vivaldi}
\address{School of Mathematical Sciences, Queen Mary,
University of London,
London E1 4NS, UK}
\email{f.vivaldi@maths.qmul.ac.uk}
%
\begin{abstract}
In rotations with a binary symbolic dynamics, a critical curve is 
the locus of parameters for which the boundaries of the partition
that defines the symbolic dynamics are connected via a prescribed 
number of iterations and symbolic itinerary.
We study the arithmetical and geometrical properties of these curves
in parameter space.
\end{abstract}
\date{\today}
\maketitle

\vspace{-20pt}
\begin{center}
{\em Dedicated to the memory of Uwe Grimm.}
\end{center}

\section{Introduction}\label{section:Introduction}

We consider a rotation $\rg$ on the circle 
(unit interval): 
\begin{equation}\label{eq:Map}
\rg:[0,1)\to[0,1)\qquad x\mapsto \{x+\theta\}, 
\end{equation}
where $\theta$ is the angle of rotation and $\{\cdot\}$ denotes 
the fractional part.
The partition of the circle
\begin{equation}\label{eq:Iab}
\gI_a=[0,\rho)\qquad \gI_b=[\rho,1),
\end{equation} 
defines a symbolic dynamics in two letters $a$ and $b$.
(For background on symbolic dynamics, see \cite{PytheasFogg, Berthe,AlessandriBerthe}).

The parameter space of this system is the closed unit square $[0,1]^2$ of 
all pairs $\zeta=(\theta,\rho)$, and a \textbf{critical point} is a pair 
$\zeta$ for which the equation
\begin{equation}\label{eq:CriticalPoint}
\ri\theta=\rj+\rho
\end{equation}
has a solution $\rz=\rij\in\Z^2$.
This means that at a critical point there is a \textbf{critical orbit} of $\rg$ containing 
both boundary points $0$ and $\rho$ of the partition (\ref{eq:Iab}). 
The \textit{shortest} portion of orbit connecting such boundary points, 
in some order, will be called the \textbf{centre} of the orbit, with the
convention that it
includes the initial boundary point but not the final one. 
We then consider the symbol sequence $w$ of the centre.
If $\theta$ is irrational, then the boundary points are visited 
only once and in a defined order, and $w$ is determined by $\zeta$. 
If $\theta$ is rational, then the orbit is periodic, and we may choose
either $0$ or $\rho$ as the initial point of the centre. 
Hence when $\theta$ is rational, $\zeta$ determines two centres and their corresponding words $w$, 
typically of different length. 
A \textbf{critical word} $w$ is a word constructed in this fashion,
and the \textbf{critical curve} $\cC_w$ is the set 
of critical points which share the same critical word $w$. 

Even though critical curves of rotations are fairly basic objects, 
to the best of our knowledge they have not been considered 
explicitly.
It is known that critical points determine the complexity of rotational 
words ---the symbol sequences generated by the rotation (\ref{eq:Map}) with
partition (\ref{eq:Iab}).
At a critical point with irrational $\theta$, a critical word of 
length $n$ is a factor (finite sub-word) of an infinite word with 
complexity $\mathcal{K}(t)$ (number of factors of length $t$) equal to 
\begin{equation}\label{eq:K(t)}
\mathcal{K}(t)=\begin{cases} 2t& t\leqslant n \\t+n & t>n,\end{cases}
\end{equation}
see \cite[theorem 10]{AlessandriBerthe}.
The case $n=1$ is the much studied \textrm{Sturmian} words 
\cite[section 6]{PytheasFogg}, while $n=\infty$ corresponds to
non-critical parameters ---the boundary points are not on the same orbit.
Note that, due to the minimality of irrational rotations, the 
collections of all factors of an infinite word (the symbolic language) 
does not depend on the initial condition \cite[p.~105]{PytheasFogg}.
An irrational critical point can therefore be characterised in 
terms of the complexity of the word of any orbit.

Critical curves are relevant to piecewise-defined dynamical
systems supporting rotational motions, such as the much-studied family 
of planar maps 
\begin{equation}\label{eq:F}
\rF(x,y)=\begin{cases} 
(\a x-y,x)& \,\,x\geqslant 0\\
(\b x-y,x)& \,\,\mbox{otherwise,}
\end{cases}
\end{equation}
where $\a$ and $\b$ are real parameters
\cite{Devaney,BeardonBullettRippon,LagariasRainsI,
LagariasRainsII,LagariasRainsIII,Garcia-MoratoEtAl}.
Since $\rF$ sends rays to the origins to themselves, 
part of its dynamics is described by a circle map with 
well-defined rotation number $\theta(\a,\b)$ which, when 
irrational, allows a topological conjugacy between the 
circle map and the rigid rotation (\ref{eq:Map}).
Each orbit of $\rF$ applied to a ray from the origin then 
corresponds to acting on the ray a growing product of 
the two matrices $A$ and $B$ defined by (\ref{eq:F}).
This product corresponds to a growing word in the symbols $A,B$,
which follows the rotational word in $\a,\b$ from the 
induced circle map.  
In this way, $F$ can be interpreted \cite{LagariasRainsI,
LagariasRainsII,LagariasRainsIII}
as the discrete 
Schr\"odinger equation on the one-dimensional 
lattice with a two-valued potential sequence generated by a
rotational word. 
Such models with quasiperiodic two-letter words developed 
via two-letter substitution rules have received 
much attention \cite{BaakeGrimm,BaakeGrimmJoseph}.

In \cite[theorem 2.2]{LagariasRainsII} it was shown that for any parameter
pair $(\a,\b)$ for which the partition boundaries 
(the positive and negative ordinate semi-axes) belong to the
same $\rF$-orbit, invariant curves exist, consisting of the union
of arcs of conic sections. In our terminology, these are critical 
points in parameter space, each with an associated finite symbol 
sequence. The locus of critical parameters in $(\a,\b)$ space
with the same sequence 
defines an algebraic curve, examples of which were first given
(not under the name of critical curves)
in \cite[examples 4.1--3]{LagariasRainsII}.
Critical curves have also appeared implicitly in the works 
\cite{SimpsonMeiss,Simpson:17,Simpson:18,Garcia-MoratoEtAl}
on mode-locking and bifurcations in piecewise-linear maps, 
where they are used to characterise sequences of mode-locking 
regions near a bifurcation point.


In \cite{RSV}, we investigated some properties of the critical curves 
of the map $\rF$, and we found structures of considerable complexity.
We believe that understanding these curves in the simpler 
---yet non-trivial--- case of 
rotations is an essential pre-requisite for the development of 
a theory of these dynamical objects in a general setting. 
This is the motivation for the present paper.

We summarise the contents and main results of this paper.
In the next section we establish terminology and notation,
and introduce the notion of a \textbf{chain}, the set of 
solutions of equation (\ref{eq:CriticalPoint}) for fixed $\rij$.
In section \ref{section:WordsOnChains} we show that 
a chain is partitioned
by Farey fractions (see \cite[chapter III]{HardyWright}) into 
a sequence of critical curves (open segments) 
separated by degenerate curves called \textbf{Farey points},
and that along a chain the code changes from $a^n$ to $b^n$, 
or vice-versa (theorem \ref{thm:CriticalCurves}). 
Furthermore, the identification of a Farey point on a chain is a relative concept:
every such point belongs to a critical curve transversal to the given chain 
(corollary \ref{cor:Endpoints}).

In section \ref{section:RationalPoints} we describe all curves through
a critical rational point, for which equation (\ref{eq:CriticalPoint}) has 
infinitely many solutions.
We show that in general every such point $\z$ is an interior point 
of exactly two curves 
---the \textbf{dominant curves} of $\z$--- as well as the common
Farey point of four infinite pencils of curves 
(theorem \ref{thm:RationalPoints}).
The symbolic dynamics of these curves is determined in
theorem \ref{thm:Words}. 

In section \ref{section:TriplePoints} we consider the geometric figure
determined by the six dominant lines of a rational critical point
and two points closest to $\z$ which share the value of $\theta$.
We show that these lines form two triples concurrent to the 
two \textbf{triple points} of $\z$, whose rotation numbers are
convergents of the continued fraction expansion of $\theta$
(theorem \ref{thm:TriplePoints}). 

\noindent{\sc Acknowledgements}.
This research was supported by the Australian Research Council grant
DP180100201 and by JSPS KAKENHI Grant Numbers JP16KK0005 and JP22K12197.

\section{Basic properties}\label{section:BasicProperties}

The parameter space for the symbolic dynamics of the map $\rg$ 
of (\ref{eq:Map}) is the set of pairs $\zeta=(\theta,\rho)\in[0,1]^2$. 
From (\ref{eq:Iab}) we find that the values $\rho=0$ and $\rho=1$ 
correspond to the trivial symbolic dynamics built from the single letter 
$b$ and $a$ respectively. 
For $\theta=0$ we have the identity map, while $\theta=1$ is included 
for consistency with Farey sequences\footnote{In phase space, the
values $\theta,\rho=1$ are represented by their fractional
part $0$, according to (\ref{eq:Map}).}.

We begin by taking a closer look at equation (\ref{eq:CriticalPoint}).
For given $\zeta$, every solution $\rij$ correspond bi-uniquely to an 
orbit segment of length $|\ri|$ having $0$ and $\rho$ as end-points, 
in the order prescribed by the sign of $\ri$. 
We first deal with the associated symbolic dynamics.

\noindent \textbf{Def}. 
A \textbf{critical orbit} of (\ref{eq:Map}) is an orbit
containing both boundary points $0$ and $\rho$, and
a \textbf{boundary word} is the symbolic dynamics of a 
finite section of a critical orbit connecting such boundary 
points, in some order, in such a way that the initial point is included and the final one is not.
A \textbf{critical word} is a boundary word of minimal length, that is,
the symbolic dynamics of the \textbf{centre} of the orbit. (Thus every boundary
word contains a critical word as a prefix.) 

If $\rho=0$, then equation (\ref{eq:CriticalPoint}) has the trivial solution 
$\rz=(0,0)$ for every $\theta$, the centre of the critical orbit is 
empty, and the critical word is the empty word $w=\varepsilon$.
There are also non-trivial solutions for rational $\theta$ 
---with associated boundary words---
which we shall consider below. Likewise, for $\rho=1$ we get the trivial
solution $\rz=(0,-1)$ as well as nontrivial solutions. 

Assume for the moment that $\rho\not=0,1$.
Then at a critical point, $\theta$ and $\rho$ are both rational or
both irrational ---indeed they belong to the same number field.
Suppose that $\theta$ is irrational.
Then the points $x=0$ and $x=\rho$ appear only 
once in the doubly infinite orbit through $0$, and hence 
(\ref{eq:CriticalPoint}) has only one solution $\rz=\rij$.
We now keep this solution fixed and regard (\ref{eq:CriticalPoint}) 
as an equation for $(\theta,\rho)$, subject to the 
constraint $(\theta,\rho)\in[0,1]^2$ (the condition 
$\rho\not=0,1$ was needed to determine $\rij$ from $\zeta$, 
and is no longer required ---see below).
We obtain a line segment of critical points ---called a \textbf{chain}---
given by
\begin{equation}\label{eq:Lij}
\cLij=\{(\theta,\rho)\,:\,\rho=\ri\theta-\rj\,,\,
\theta^-\leqslant\theta\leqslant\theta^+\}
\end{equation}
where
\begin{equation}\label{eq:Thetapm}
\theta^-= \begin{cases} \rj/\ri& \ri >0\\
                     (\rj+1)/\ri& \ri<0
       \end{cases}
\qquad\mbox{and}\qquad \theta^+=\theta^-+\frac{1}{|\ri|}.
\end{equation}
Note that $\ri\not=0$ by assumption and that $(\ri,\rj)$ is also a 
solution of equation (\ref{eq:CriticalPoint}) when $\rho=0,1$,
namely at $(\theta^{-s},0)$ and $(\theta^s,1)$,
where $s=\mathrm{sign}(\ri)$.

The above association between a solution $\rij$ of (\ref{eq:CriticalPoint}) and
the chain $\cLij$ is extended to the cases $\rho=0,1$, by defining
\begin{equation}\label{eq:L0}
\cL_{0,0}=\{(\theta,0)\,:\, 0\leqslant \theta\leqslant 1\}
\qquad
\cL_{0,-1}=\{(\theta,1)\,:\, 0\leqslant \theta\leqslant 1\}.
\end{equation}
From the above and (\ref{eq:CriticalPoint}) one verifies that 
for any integer $\ri$, a pair $\rij\in\Z^2$ corresponds to a chain 
if $\rj$ is subject to the bounds 
\begin{equation}\label{eq:Jbounds}
J(\ri)=
\begin{cases} 0\leqslant \rj\leqslant\ri-1 & \ri>0\\
             -1\leqslant \rj \leqslant 0 &  \ri=0\\
              \ri\leqslant\rj\leqslant -1& \ri<0.
\end{cases}
\end{equation}
A pair of integers $\rij$ satisfying the above conditions
will be called the \textbf{affine parameters} of the chain $\cLij$.
They are the solution of (\ref{eq:CriticalPoint}) shared by
all points on the chain. 

We are interested in the critical words $w=w(\theta,\rho(\theta))$ 
of the points of $\cL_{\ri,\rj}$ of (\ref{eq:Lij}) and (\ref{eq:L0}). 
The following definition relates points and words.

\noindent \textbf{Def.} Let $\cLij$ be a chain, let $\z\in\cLij$, and
let $w$ be the critical word at $\z$ (not necessarily of length $|\ri|$).
The \textbf{critical curve $\cC_w$ on $\cLij$ containing $\z$}
is the set of points of $\cLij$ sharing the same critical word. 
If such a set reduces to a point, we shall speak of a 
\textbf{Farey point} $\cE_w$.

Thus a chain is partitioned into critical curves and Farey points.
In particular, each of the chains (\ref{eq:L0}) has empty word
as critical word, and therefore consists of a single critical 
curve $\cC_\varepsilon$, and no Farey points.

\noindent \textbf{Def.} A non-empty boundary word $w$
is \textbf{positive} (\textbf{negative}, respectively) 
if the first symbol of $w$ is $a$ ($b$, respectively).
The empty word is both positive and negative.
Likewise, the sign of the affine parameters $\rz=\rij$ is 
that of $\ri$ if $\ri\not=0$, and if $\ri=0$ then $\rz$ 
is both positive and negative.

We see that a non-empty word is positive (negative) if the corresponding 
centre starts at $0$ ($\rho$). Let $\rho\not=0,1$.
Because ---as noted above--- to an irrational critical points $\z$ 
there corresponds unique affine parameters, the point $\z$ has 
a well-defined sign and so does the unique critical curve that contains it. 
By the same criterion, rational critical points are both 
positive and negative.

Let $n=|\ri|$.
At all irrational points on $\cLij$ the critical word 
$w=w_0\cdots w_{n-1}$ has length $n$, but this is not 
necessarily the case if $\zeta$ is rational. 
Indeed in this case the orbit is periodic and, being also 
critical, both points $0$ and $\rho$ are visited infinitely often.
As a result, equation (\ref{eq:CriticalPoint}) has a 
doubly-infinite set of solutions $(\ri_t,\rj_t)$, which
include $(\ri,\rj)$, and to each solution there is an associated
boundary word. If there is a $t$ such that $\ri_t$ and $\ri$ have 
the same sign and $|\ri_t|<n$,
then the centre of the critical orbit is shorter than $n$,
and therefore $w$ is not a critical word.
In what follows, when we speak of the boundary words of a chain
$\cLij$, we will always refer to words of length $|\ri|$.

From the above discussion we conclude that a rational critical 
point should be regarded as being both positive and negative,
and this duplicity is reflected in the sign of the boundary
words at that point. Such a correspondence however fails at the
rational points of the chains (\ref{eq:L0}), where the (non-empty)
boundary words assume only one sign, since one element of the partition 
(\ref{eq:Iab}) is empty. 
Specifically, at the rational points $(p/q,0)$, equation 
(\ref{eq:CriticalPoint}) has the solutions $\rij=(tq,tp),t\in\Z$,
assuming both signs. However, for $t\not=0$ the corresponding
boundary words $b^{|tq|}$ are negative.
There is an analogous discrepancy at $(p/q,1)$, where all 
boundary words are positive.

\section{Boundary and critical words on a chain}\label{section:WordsOnChains}

\begin{figure}[tb]
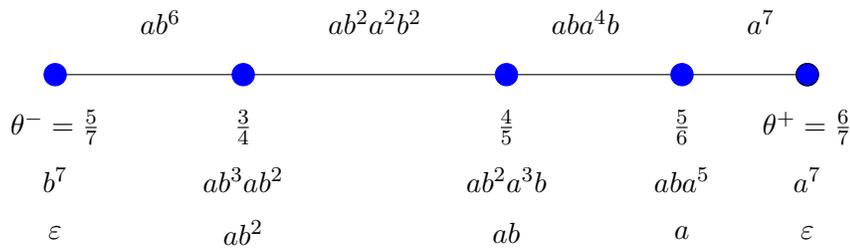

\centering
\include{fig-Chain3}
\caption{\label{fig:Chain}\rm\small
The partition of the positive chain $\cL_{7,5}$ with affine 
parameters $\rij=(7,5)$, into four critical curves and five 
Farey points (solid circles), determined according to 
theorem \ref{thm:CriticalCurves} i). 
Along the chain, all boundary words have length $\ri=7$.
Above the line we have the critical words of the curves,
and below the line the Farey fractions with corresponding 
boundary words
[see theorems \ref{thm:CriticalCurves} ii) and \ref{thm:Words}].
The bottom row displays the critical words at the Farey points, 
all of length smaller than $7$, including the empty word 
$\varepsilon$ of zero length at $\theta^\pm$.
The boundary word at a Farey point is the concatenation of 
a critical word and a periodic word, whose period is given by
the denominator of the fraction.
}
\end{figure}

In this section we consider the decomposition of a chain into critical 
curves and Farey points, as defined in section \ref{section:BasicProperties},
as well as the associated symbolic dynamics.
Figure \ref{fig:Chain} serves as an illustration of the items we shall
be dealing with.

Our first result describes all the boundary words on a chain.

\begin{theorem}\label{thm:CriticalCurves}
Let $\cLij$ be a chain, let $n=|\ri|\geqslant 1$ and let
$\cF_n$ be the $n$th Farey sequence. Then

\noindent i)\qad 
The set of rotation numbers
\begin{equation}\label{eq:gFarey}
\cF=\cF_n\cap [\theta^-,\theta^+]
\end{equation}
partitions $\cLij$ 
into $|\cF|-1$ critical curves, themselves segments,
separated by $|\cF|$ Farey points.

\noindent ii)\qad The boundary words on $\cLij$ are computed
recursively as follows. Assume first that $\cLij$ is positive,
and let $w=w_0\cdots w_{n-1}$ be the positive word of length 
$n$ at $\theta=p/q\in\cF$. Finally, let $w^\pm$ 
be the words of the adjacent curves to the right ($+$) and
left ($-$) ($w^\pm$ is missing 
at $\theta=\theta^\pm$). 
Then, for $k=1,\ldots,n-1$ the following holds:
\begin{enumerate}
\item [1)] At $\theta=\theta^-$ we have $w=b^n$, and
$w_k^+=a$ iff $k\equiv 0\mod{q}$. This holds also for $k=0$.
\item [2)] At $\theta=\theta^+$ we have $w=a^n$, and
$w_k^-=b$ iff $k\equiv 0\mod{q}$.
\item [3)] For $\theta\in\cF\setminus \{\theta^-,\theta^+\}$ we have
\begin{enumerate}
\item [$\alpha$)] $w_k=b$ and $w_k^-=w_k\not =w_k^+$ iff $k\equiv n\mod{q}$;
\item [$\beta$)] $w_k=a$ and $w_k^-\not=w_k=w_k^+$ iff $k\equiv 0\mod{q}$.
\end{enumerate}
\end{enumerate}
The corresponding statements for negative chains are obtained 
from the above by exchanging all $a$s and $b$s.
\end{theorem}

The statement of the theorem excludes the case $\ri=0$. We remark
that the critical word for the chains (\ref{eq:L0}) is $\varepsilon$
by definition, while the boundary words at the rational points are 
given in parts ii) 1,2). 

\proof
i) Along the chain (\ref{eq:Lij}), the point $x_t(\theta)$,
$t=0,\ldots,n$ of the centre are affine functions of $\theta$
with slope $t$ if $w$ is positive, and slope $-n+t$ if $w$ is negative.
If for some $\theta$ in the range (\ref{eq:Thetapm}) two such functions 
coincide, then the map $\rg$ is periodic with period not exceeding $n$,
that is, $\theta\in\cF$. Conversely, any $\theta\in\cF$
corresponds to a periodic orbit of $\rg$ of period not exceeding $n$.
Note that $\theta^\pm$ are the only elements of $\mathcal{F}$ whose 
denominator is divisible by $n$, because the numerators of $\theta^\pm$ 
are consecutive integers. 

Assume first that $\theta=p/q\in\cF\setminus \{\theta^-,\theta^+\}$. Then $1<q<n$,
because $q$ does not divide $n$.
Periodicity implies that $x_0(p/q)=x_q(p/q)$, and the intersection
of $x_0$ and $x_q$ is transversal because the slopes of the two 
functions are different. It follows that the $q$th symbol of the critical word
$w$ changes at $p/q$, which is therefore the common end-point of two adjacent
segments, hence a Farey point. 
Since $\theta^\pm$ are necessarily Farey points, and all rationals
with denominator not exceeding $n$ have been accounted for, we have shown 
that the elements of $\cF$ partition the chain $\cLij$
into $|\cF|-1$ curves, as desired.

ii) Let $w$ be positive.

1) At $\theta=\theta^-$ we have $\rho=0$; then, according to (\ref{eq:Iab}),
the partition element $\gI_a$ is empty, so the code is $b^n$.
We have $x_k=0$ iff $k\equiv 0\mod{q}$,
and these are precisely the values of $k$ (which include $k=0$) 
for which $x_k\in\gI_a$ for $\theta>\theta^-$, that is, $w_k^+=a$.

2) At $\theta=\theta^+$ we have $\rho=1$; then $\gI_b$ is empty, 
and the code at $\theta^+$ is $a^n$. 
Again we have $x_k=0$ (on the circle) iff $k\equiv 0\mod{q}$, and if $k\not=0$
we have $x_k\in\gI_b$ for $\theta<\theta^+$, that is, $w_k^+=b$.
(The value $k=0$ must be excluded here, because $w_0=a$ for
all $\theta\not=\theta^-$.)

3) Let now $\theta=p/q\in\cF\setminus \{\theta^-,\theta^+\}$.
Then $n\not\equiv0\mod{q}$, as noted above.

$\alpha$) \/ The critical curve property together with $q$-periodicity implies 
that $x_{k}=x_n=\rho$ iff $k\equiv n\mod{q}$. 
The proper symbol $w_k$ at $\theta$ is $b$, and since $k<n$,
the slope of $x_k$ is smaller than that of $x_n$, so that
$x_k\in\gI_b$ ($x_k\in\gI_a$) in a left (right) 
neighbourhood of $\theta$, that is, $w_k^-=w_k\not=w_k^+$. 

$\beta$) \/ $q$-periodicity implies that $x_k=x_0=0$ iff $k\equiv 0\mod{q}$. 
The proper symbol $w_k$ at $\theta$ is $a$, and since $k>0$,
the slope of $x_k$ is greater than that of $x_0$, so that
$x_k\in\gI_b$ ($x_k\in\gI_a$) in a left (right) 
neighbourhood of $\theta$, that is, $w_k^-\not=w_k=w_k^+$. 

The proof of ii) is complete.

If $w$ is negative, the argument develops in a symmetrical manner.
At $\theta^+$ we have $\rho=0$ whence $w=b^n$.
As $\theta$ decreases, all collisions of orbit points involve
$b$s turning into $a$s, until we reach $\theta^-$ with code $a^n$.
We omit the details.
\endproof


The arithmetical and combinatorial aspects of a chain
are illustrated in figure \ref{fig:Chain}.
Applying theorem \ref{thm:CriticalCurves} ii) recursively along 
a chain, from $\theta^-$ to $\theta^+$, we deduce that 
every letter $b$ of the initial word changes to an $a$ without 
omissions or repetitions. 
This translates into the following arithmetical statement.

\begin{corollary}\label{cor:Congruences}
For any integers $n>m$ and $m \geqslant 0$, let $\mathcal{F}$ be 
the subset of $\cF_n$ lying between $m/n$ and $(m+1)/n$, and for 
each $p/q\in\mathcal{F}$ consider the congruences
$$
x\equiv n\mod{q},\qquad x\equiv 0\mod{q}
$$
(which coincide if $q$ divides $n$).
Then, as $p/q$ ranges in $\mathcal{F}$, the solutions of this family 
of congruences form a complete set of residues modulo $n$, and 
each non-zero residue is a solution of exactly one congruence.
\end{corollary}
In the above statement there is no restriction on the numerator $m$, 
because the restriction that appears in theorem \ref{thm:CriticalCurves} 
plays no role in the proof of part ii).

The number $M=|\mathcal{F}|-1$ of curves in a chain depends on its 
affine parameters $(\ri,\rj)$.
Let $n=|\ri|$. Such a number is independent of $n$ in only
three cases, namely $\rj=0$, $\rj=n-1$ ($M=1$), and $n \geqslant 3$ odd 
and $2\rj+1=n$ ($M=2$), plus the corresponding values for negative $\ri$. 
In all other cases, for fixed $\rj$,
we have  $M(\ri,\rj)\to\infty$ as $|\ri|\to\infty$. In this case the 
average order of $M(n)$ is $3n/\pi^2$, which may be deduced from 
that of the Farey series \cite[theorem 331]{HardyWright}.

Our next result provides alternative characterisations of 
the Farey points of a chain.

\begin{lemma}\label{lma:Endpoints}
Let $\cLij$ be a chain and let $\z=(\theta,\rho)\in\cLij$.
The following statements are equivalent:
\vspace*{-10pt}
\begin{enumerate}
\item [i)]   $\z$ is a Farey point of $\cLij$;
\item [ii)]  the critical word $w$ at $\z$ is such that $|w|<|\ri|$;
\item [iii)] there exists $(\ri',\rj')$ with $|\ri'|<|\ri|$ and
$\mathrm{sign}(\ri')=\mathrm{sign}(\ri)$ such that
$\z=\cLij\cap\cL_{\ri',\rj'}$
and $\z$ is not a Farey point of $\cL_{\ri',\rj'}$. 
\end{enumerate}
\end{lemma}
\proof 
For $\ri=0$ all statements above are false, hence equivalent. 
If $\ri\not=0$ and $\z$ also belongs to $\cL_{0,0}$ or $\cL_{0,-1}$,
then all statements are true since the critical word of these
chains is the empty word, which is both positive and negative by definition.

We now assume that $\ri\not=0$, $\rho\not=0,1$, and we
let $w$ be the critical word at $\z$. 
We shall prove that i) $\Rightarrow$ iii) $\Rightarrow$ ii) $\Rightarrow$ i).

i) $\Rightarrow$  iii).
If $\z$ is a Farey point of $\cLij$, then from theorem 
\ref{thm:CriticalCurves} i) and the fact that 
$\rho\not=0,1$ we have that $\theta=p/q$ with $q<|\ri|$ and $q \nmid |\ri|$.
We define $c:= \lfloor |\ri| / q \rfloor \geqslant 1$ and $\ri'$ and $\rj'$ by $\ri'=\ri-\mathrm{sign} (\ri)\,c\,q$ and  
$\rj'=\rj-\mathrm{sign}(\ri)\,c\,p$. Hence $1 \leqslant |\ri'| < q < |\ri|$, $\mathrm{sign}(\ri')=\mathrm{sign}(\ri)$ and one checks that $\ri'\,\theta - \rj'=\ri\,\theta - \rj=\rho$.
So $\z$ lies at the 
intersection of $\cLij$ and $\cL_{\ri',\rj'}$. Furthermore, if $w$ is the critical word at $\z$
with the same sign as $\ri$, it has the minimal
length $|\ri'|$ by construction. Since $q > |w|=|\ri'|$, theorem \ref{thm:CriticalCurves} i) shows that $\z$
cannot be a Farey point on $\cL_{\ri',\rj'}$.

iii) $\Rightarrow$  ii).
If iii) holds, then $\theta=(\rj-\rj')/(\ri-\ri')$,
and hence at $\zeta$ the critical orbit is periodic with 
period $|\ri|-|\ri'|$. Since the length of the critical word
is necessarily smaller than the period, we have
$|w|<|\ri|-|\ri'|<|\ri|$. 

ii) $\Rightarrow$ i). 
If ii) holds, then $\z$ is necessarily
a rational point, and the denominator of $\theta$ is less than
$|\ri|$, being a divisor of $|\ri|-|w|$. Thus $\theta\in\cF_{|\ri|}$,
and $\z$ is a Farey point from theorem \ref{thm:CriticalCurves} i).
\endproof

We infer from part iii) of the lemma, 
which itself relies on theorem \ref{thm:CriticalCurves}, that Farey points exist only in the context of a given chain. More precisely, at a Farey point $\cB_w$ on a chain
there is always a direction in parameter space along another chain for which the  
critical word $w$ does not change for sufficiently small displacements (and so 
$\cB_w$ is not a Farey point on the second chain).
This fact is expressed concisely by the following statement:
\begin{corollary}\label{cor:Endpoints}
A Farey point of a chain belongs to a critical curve transversal to
the chain.
\end{corollary}

\section{Curves through rational points}\label{section:RationalPoints}

We now characterise all critical curves through a rational critical 
point $\zeta$. As discussed in the previous section, equation
(\ref{eq:CriticalPoint}) for rational $\zeta$ has a doubly-infinite family
of solutions corresponding to as many chains passing through $\zeta$.
We shall identify the chains for which $\zeta$ belongs to a critical curve,
and those for which $\z$ is a Farey point; in the latter case we also 
determine the adjacent curves on the chain. 
In what follows, by the \textbf{Farey points of a critical curve} we shall 
mean those belonging to the closure of the curve. 

If $\theta=p/q$, then the orbit of $0$ consists of
$q$ equally spaced points on the unit interval, and therefore
$\z$ must be of the form $(p/q,r/s)$, with $s$ dividing $q$.
We call $q$ the \textbf{denominator} of $\z$.
We develop $p/q$ in continued fractions, and of
the two possible continued fractions representations 
we choose the one whose last coefficient is unity 
(see \cite[theorem 162]{HardyWright}). Then the index $n$
of the last convergent $p_n/q_n=p/q$ is determined unambiguously.
At various junctures we shall discuss the implications of this
choice.
We let 
\begin{equation}\label{eq:tau}
\tau=q'\frac{r}{s},\qquad 
   \tau^\pm =q'\biggl(\frac{r}{s}\pm\frac{1}{q}\biggr)
\qquad\mbox{where}\qquad q'=(-1)^{n-1}q_{n-1}.
\end{equation}

Given a rational point $\z$ with $\rho\not=0,1$, we consider,
among the positive chains containing $\z$ that with 
minimal value of $\ri$, and denote its affine parameters 
by $(\bi^+,\bj^+)$. In a similar manner, the negative chain with
minimal value of $-\ri$ will be denoted by $(\bi^-,\bj^-)$.
These pairs will be called the \textbf{dominant affine parameters} 
of the point $\z$. Likewise, we shall speak of the
\textbf{dominant chains} (or \textbf{lines}),
and the \textbf{dominant curves} of $\z$. 
We will show that these objects exist and are unique. 
For uniformity of exposition, we treat the cases $\rho=0,1$ as
follows. We shall regard the segment $\rho=0$ as the positive dominant
chain of any point $\zeta=(\theta,0)$, and the segment $\rho=1$ as the
negative dominant chain of any point $\zeta=(\theta,1)$. In both
cases, the chain consists of a single curve with the empty
word.
Accordingly, if $\rho=0$, we let $(\bi^+,\bj^+)=(0,0)$ and
$(\bi^-,\bj^-)=(-q,-p)$.
If $\rho=1$ we let $(\bi^+,\bj^+)=(q,p-1)$ and $(\bi^-,\bj^-)=(0,-1)$.
Finally, we define the \textbf{upper} and \textbf{lower neighbours} 
of a rational point $\z=(p/q,r/s)$ to be the points 
\begin{equation}\label{eq:zud}
\zu=\z+\Bigl(0,\frac{1}{q}\Bigr)
    =\Bigl(\frac{p}{q},\frac{r}{s}+\frac{1}{q}\Bigr)
\qquad
\zd=\z-\Bigl(0,\frac{1}{q}\Bigr)
    =\Bigl(\frac{p}{q},\frac{r}{s}-\frac{1}{q}\Bigr),
\end{equation} 
respectively. 
If $\z$ is a critical point, then so are $\zu$ and $\zd$, 
with one of them missing if $\rho=0,1$.

We shall consider the four 
quadrants with origin at $\z$ and label all objects that are 
pertinent to such quadrants ---points and curves--- with the 
superscripts $\I,\II,\III,\IV$.
The superscript $\pm$ will refer to the sign of curves, so that
the positive sign refers to $\I,\III$ and the negative sign 
to $\IV,\II$. 

\begin{figure}[t]
\centering
         \hspace*{70pt}\includegraphics[scale=0.80,clip]{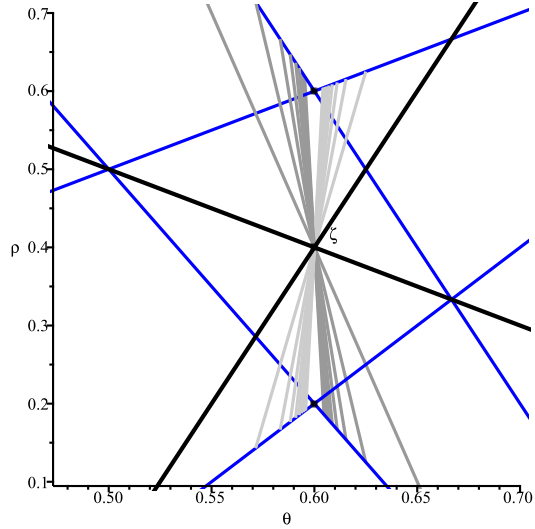}
\vspace*{-400pt}
\caption{\label{fig:AllArcs} 
The critical curves through (or adjacent to) the rational point 
$\z=(3/5,2/5)$ (see theorem \ref{thm:RationalPoints}), with its two
dominant lines (black) and four pencils of curves concurrent
at $\z$ (grey), the latter being their common Farey point.
The other Farey point of the curves in a pencil lie on a dominant
lines of one of $\z$'s neighbours (blue), each pencil paired
with a different line.
The six dominant lines feature two concurrent triples, as 
detailed in theorem \ref{thm:TriplePoints}.
}
\end{figure}

\goodbreak
\begin{theorem}\label{thm:RationalPoints} 
Let $\z=(\theta,\rho)$ be a rational critical point.
If $\rho\not=0,1$ then $\z$ is an interior point of the dominant
curves of $\z$, and the common Farey point of four infinite 
pencils of curves, arranged pairwise as adjacent curves in two
pencils of chains of opposite sign.
Furthermore, all Farey points distinct from $\z$ in a pencil belong 
to (the closure of) a dominant curve
of a neighbour of $\z$, this association being
bi-unique ---see figure \ref{fig:AllArcs}.
The same applies if $\rho=0,1$, but in this case one neighbour
and two pencils are missing (three pencils, if $\theta=0, 1$).
\end{theorem}

\proof 
By assumption, $\z$ has the form $(p/q,r/s)$ with $s$ dividing $q$
and it belongs to some curve or Farey point;
we begin to construct all chains through $\z$.
From (\ref{eq:CriticalPoint}) 
with $q=su$ we find $ru=\ri p -\rj q$, with solutions
\begin{equation}\label{eq:Solutions}
\ri_t=ru q'+tq,\quad \rj_t=ru p'+tp\qquad t\in\Z,
\end{equation}
where
$p'=p_{n-1}(-1)^{n-1}$ and $q'=q_{n-1}(-1)^{n-1}$.
Thus $pq'-qp'=1$. 

The sequences (\ref{eq:Solutions}) are the affine parameters of
all chains containing $\z$.
Assume first that $\rho\not=0,1$, that is, $s\not=1$. 
If we had $\ri_t=0$ for some $t$,
then $s$ would also divide $q_{n-1}$, which in turn would yield $s=1$,
contrary to the assumption. 
So there are exactly two values $t^\pm$ of $t$ in (\ref{eq:Solutions}) 
for which $0<|\ri_t|< q$, namely
\begin{equation}\label{eq:tpm}
t^+=\lceil -q'r/s\rceil,
\qquad
t^-=\lfloor -q'r/s\rfloor=t^+-1
\end{equation}
while all other values of $t$ give $|\ri_t|>q$.

With reference to (\ref{eq:tpm}) define
\begin{equation}\label{eq:ipmfirst}
\ri^+(\ell)=\ri_{t^++\ell}
\qquad
\ri^-(\ell)=\ri_{t^--\ell}
\qquad \ell=0,1,2,\ldots
\end{equation}
and similarly for $\rj^\pm(\ell)$, so that the sign of 
$\ri^\pm(\ell)$ is $\pm1$ for all $\ell$. 
To obtain explicit formulae for $\ri^\pm(\ell)$ and $\rj^\pm(\ell)$
as functions of $p/q$ and $r/s$, we use 
(\ref{eq:Solutions}--\ref{eq:ipmfirst}), 
keeping in mind that $ru=qr/s$ and $\lceil -x\rceil=-\lfloor x\rfloor$. 
We obtain the dominant affine parameter of $\z$:
\begin{equation}\label{eq:DominantIJ}
\bi^\pm=\ri^\pm(0)=\pm q\{\pm \tau\},
\qquad
\bj^\pm=\rj^\pm(0)=\pm p\{\pm \tau\}-\frac{r}{s},
\end{equation}
where $\{\cdot\}$ denotes the fractional part, and $\tau$ was
defined in (\ref{eq:tau}).
Finally, we obtain, for $\ell=0,1,2,\ldots$
\begin{equation}\label{eq:ijpm}
\ri^\pm(\ell)=\pm q(\{\pm \tau\}+\ell)
\qquad
\rj^\pm(\ell)=\pm p(\{\pm\tau\}+\ell)-\frac{r}{s}.
\end{equation}
The above expression give all affine parameters of all chains 
though $\z$ in the case $\rho\not=0,1$.
From (\ref{eq:Lij}) the sign of $\ri_t$ is the sign of 
the chain containing $\z$, and $|\ri_t|=|w|$. Then, formulae
(\ref{eq:ijpm}) and lemma \ref{lma:Endpoints} give two 
infinite sequences of chains of opposite sign, such that 
$\z$ is an interior point of 
a curve in the chains $[\ri^\pm(0),\rj^\pm(0)]$, and a Farey point in 
all other chains, as claimed.

We now assume further that $\rho\not=1/q, (q-1)/q$.
As indicated above, we shall use the superscripts $\I,\II,\III,\IV$, 
which refer to the four quadrants with origin at $\z$, to label the 
points in the four sequences and other relevant quantities.
We denote the dominant affine parameters of the upper neighbour
$\zu$ ---see (\ref{eq:zud})--- by
$(\bi^\I,\bj^\I)$ and $(\bi^\II,\bj^\II)$, and those of $\zd$ by
$(\bi^\III,\bj^\III)$ and $(\bi^\IV,\bj^\IV)$. (Thus $\I,\III$ are
positive and $\II,\IV$ are negative.)
Considering (\ref{eq:tau}) and (\ref{eq:DominantIJ}), we find
\begin{equation}\label{eq:isigma}
\bi^\I=q\{\tau^+\}\qquad
\bi^\II=-q\{-\tau^+\}\qquad
\bi^\III=q\{\tau^-\}\qquad
\bi^\IV=-q\{-\tau^-\},
\end{equation}
and
$
\bj^\sigma=\bi^\sigma\frac{p}{q} -\rho^\sigma,
$
where $\rho^\sigma=r/s+1/q$ if $\sigma=\I,\II$ 
and $\rho^\sigma=r/s-1/q$ if $\sigma=\III,\IV$.

Next, for each pencil, we determine the external Farey point 
of each curve having $\z$ as the common Farey point. 
The corresponding parameters are given in (\ref{eq:ijpm}) 
with $\ell\geqslant 1$. 
Then $|\ri^\pm(\ell)|>q>|\bi^\sigma|$ for every $\sigma$, and
therefore, from lemma \ref{lma:Endpoints} iii) (with 
$(\ri',\rj')=(\bi^\sigma,\bj^\sigma)$ and $(\ri,\rj)=(\ri^\pm(\ell),\rj^\pm(\ell))$), we see that
the curves of the pencil in sector $\sigma$ cannot extend further 
than the dominant line $(\bi^\sigma,\bj^\sigma)$.
We will now show that all Farey points in fact belong to that line.

Let us consider the Farey sequence $\mathcal{F}$, given in (\ref{eq:gFarey}),
for the affine parameters $(\ri^\pm(\ell),\rj^\pm(\ell))$. 
From theorem \ref{thm:CriticalCurves} i) and the fact that
$p/q\not=\theta^\pm$, the rotation numbers of the Farey points of 
the curves adjacent to $\z$ are three successive terms of $\mathcal{F}$:
$p^l/q^l<p/q<p^r/q^r$.
Then (see \cite[section 3.1]{HardyWright}), we have $p^r q-q^r p=1$.
Moreover $p/q$ is the mediant of $p^r/q^r$ and $p^l/q^l$, and we shall 
compute the latter from the former.

With the notation as above, we have the candidate values for $p^r$ and $q^r$:
\begin{equation}\label{eq:Solutions2}
q^r_k=-q'+kq\qquad 
p^r_k=-p'+kp\qquad k\in\Z.
\end{equation}
As $|k|$ increases, $p^r_k/q^r_k$ approaches $p/q$, and 
the approach is from the right if $q^r_k$ is positive. 
The fraction $p^r_k/q^r_k$ belongs to the $|\ri^\pm(\ell)|$th Farey 
sequence if $0< q^r_k\leqslant |\ri^\pm(\ell)|$, and so we let $k^\pm$ 
be the largest value of $k$ for which this property holds.
We find
\begin{equation}\label{eq:k^+}
k^\pm(\ell)=\ell\pm t^\pm +\lfloor \pm\tau^\pm \rfloor,
\end{equation}
where $\tau^\pm$ was defined in (\ref{eq:tau}).
Using the quadrant superscripts, we let 
$$
p^\I(\ell)=p^r_{k^+},\quad
q^\I(\ell)=q^r_{k^+},\quad
p^\IV(\ell)=p^r_{k^-},\quad
q^\IV(\ell)=q^r_{k^-}. 
$$
To harmonise the notation, we shall also use the symbols
$$
\ri^\sigma(\ell)=\begin{cases} \ri^+(\ell)& \sigma=\I,\III\\
                        \ri^-(\ell)& \sigma=\II,\IV
          \end{cases}
$$
where $\ri^\pm(\ell)$ was defined in (\ref{eq:ijpm}).

To compute $p^{\II,\III}(\ell)$ and $q^{\II,\III}(\ell)$ using the mediant 
property we let $p^{\II,\III}(\ell)=a^\pm p-p^{\IV,\I}(\ell)$ and 
$q^{\II,\III}(\ell)=a^\pm q-q^{\IV,\I}(\ell)$, 
where $a^\pm$ is a positive integer. The required value
of $a^\pm$ is the largest such that $q^{\III,\II}\leqslant |\ri^\pm(\ell)|$, 
which is 
$$
a^\pm(\ell)=2(\ell\pm t^\pm)+\lfloor\pm\tau^\pm\rfloor
        +\lfloor\mp\tau^\mp\rfloor.
$$
Performing the calculation explicitly gives:
\begin{equation}\label{eq:pq}
\begin{array}{lrlclcl}
p^\I(\ell)&=&  -p'+p(\ell+t^++\lfloor\tau^+\rfloor)&\qquad&
q^\I(\ell)&=&  -q'+q(\ell+t^++\lfloor\tau^+\rfloor)\\
p^\II(\ell)&=& +p'+p(\ell-t^-+\lfloor-\tau^+\rfloor)&\qquad&
q^\II(\ell)&=& +q'+q(\ell-t^-+\lfloor-\tau^+\rfloor)\\
p^\III(\ell)&=&+p'+p(\ell+t^++\lfloor\tau^-\rfloor)&\qquad&
q^\III(\ell)&=&+q'+q(\ell+t^++\lfloor\tau^-\rfloor)\\
p^\IV(\ell)&=& -p'+p(\ell-t^-+\lfloor-\tau^-\rfloor)&\qquad&
q^\IV(\ell)&=& -q'+q(\ell-t^-+\lfloor-\tau^-\rfloor).
\end{array}
\end{equation}

We have constructed four infinite sequences of Farey points, and we must
now show that the Farey points of each sequence belong to (the closure of)
a dominant curve of a neighbour of $\z$.
We begin to show that
these points are collinear, and to this end, we let
\begin{equation}\label{eq:ell}
\theta(\ell)=\frac{p(\ell)}{q(\ell)}
\qquad
\rho(\ell)=\ri(\ell)\theta(\ell)-\rj(\ell)
\qquad
\z(\ell)=(\theta(\ell),\rho(\ell)),
\end{equation}
where all quantities refer to the same superscript.
Using the above and (\ref{eq:pq}), we find
\begin{equation}\label{eq:limits}
\lim_{\ell\to\infty}\z^\sigma(\ell)
=
\begin{cases}
\zu & \sigma=\I,\II\\
\zd & \sigma=\III,\IV
\end{cases}
\end{equation}
where $\zud$ was defined in (\ref{eq:zud}).

From (\ref{eq:ell}) and (\ref{eq:pq}) we find, 
after some manipulations,
\begin{equation}
\alpha^\sigma:=\frac{\rho(\ell+1)-\rho(\ell)}{\theta(\ell+1)-\theta(\ell)}
=\ri^\sigma(\ell)-\gamma q^\sigma(\ell) 
\qquad \gamma=\begin{cases} +1& \sigma=\I,\III\\ -1&\sigma=\II,\IV\end{cases}
\end{equation}
which is valid for $\ell\geqslant 1$.
Evaluating the last expression gives $\alpha^\sigma=\bi^\sigma$,
as in (\ref{eq:isigma}), independent of $\ell$. 
Thus the line to which the points $\z^\sigma(\ell)$ as well as the limit point belong
has affine parameters $(\bi^\sigma,\bj^\sigma)$, where
$
\bj^\sigma=\bi^\sigma\frac{p}{q} -\rho^\sigma(\infty).
$
Comparing (\ref{eq:isigma}) with (\ref{eq:ijpm}), we see that
$\sigma=\I,\II$ correspond to the dominant affine parameters 
of $\zu$, while $\sigma=\III,\IV$ are those of $\zd$.

We have proved that the Farey point distinct from $\z$ in the
$\sigma$-pencil, namely $\z^\sigma(\ell),\ell=1,2\ldots$ belong 
to the dominant line with parameters $(\bi^\sigma,\bj^\sigma)$.
It remains to show that all these Farey points belong to the closure of
the dominant curve (infinitely many of them do, since $\zud$ is 
an interior point of the dominant curve).
By theorem \ref{thm:CriticalCurves} i), we must show 
that the element of the $|\bi^\sigma|$th Farey sequence
which lies to the right (for $\sigma=\I,\IV$) or to the left
(for $\sigma=\II,\III$) of $p/q$ lies at least as far out as
$p^\sigma(1)/q^\sigma(1)$.
From (\ref{eq:ijpm}) and (\ref{eq:isigma}) we have 
$|\bi^\sigma|\leqslant q\leqslant |\ri^\sigma(1)|$, and 
hence the corresponding Farey sequences satisfy
$\mathcal{F}_{|\bi^\sigma|}
\subset 
\mathcal{F}_q
\subset
\mathcal{F}_{|\ri^\sigma(1)|}$
which proves our assertion.

We now discuss the cases $\rho=1/q, (q-1)/q$.
Let $\rho=1/q$, i.e., $\z=(\theta,\rho)=(p/q,1/q)$.
The Farey points distinct from $\z$ in the $\I$- and $\II$-pencils can
be treated in the same way as before.
Thus, we focus on those in the $\III$- and $\IV$-pencils.
Since $\rho\not=1$, we have $q \geqslant 2$, which gives $\theta\not =
0,1$, i.e., $1 \leqslant p \leqslant q-1$.
Thus, $\z$ is in the triangular region specified by $\theta - \rho
\geqslant 0$, $\theta + \rho \leqslant 1$, and $\rho > 0$.
Since $|\ri^\pm(\ell)|>q$ for $\ell\geqslant 1$, the line with
parameters $(\ri^\pm(\ell),\rj^\pm(\ell))$, $\ell\geqslant 1$ intersects
the segment $\rho=0$ ($0 \leqslant \theta \leqslant 1$), i.e., the
positive dominant curve of $\zd = (p/q,0)$.
We will now show that all the intersection points are the Farey points
distinct from $\z$ in the $\III$- and $\IV$-pencils.
The $\theta$-coordinate of the intersection point of the lines $\rho=
\ri^\pm(\ell) \theta- \rj^\pm(\ell)$ and $\rho=0$ is given by $\theta =
\rj^\pm(\ell)/\ri^\pm(\ell)$.
Using the condition for neighbouring terms in a Farey sequence, we see
that $\rj^+(\ell)/\ri^+(\ell)$ and $p/q$ (resp. $p/q$ and
$-\rj^-(\ell)/-\ri^-(\ell)$) are neighbours in
$\mathcal{F}_{|\ri^+(\ell)|}$ (resp. $\mathcal{F}_{|\ri^-(\ell)|}$).
This proves the assertion for the case $\rho=1/q$.
In the case of $\rho=(q-1)/q$, i.e., $\z=(\theta,\rho)=(p/q,(q-1)/q)$,
we can similarly show that the line with parameters
$(\ri^\pm(\ell),\rj^\pm(\ell))$, $\ell\geqslant 1$ intersects the
segment $\rho=1$ ($0 \leqslant \theta \leqslant 1$), i.e., the
negative dominant curve of $\zu = (p/q,1)$ and that all the intersection
points are the Farey points distinct from $\z$ in the $\I$- and
$\II$-pencils.

The final item in the proof are the boundary cases $\rho=0,1$.
By definition, a rational critical point $\z$ with $\rho=0$ has the
following affine parameters for positive and negative chains: for
$\ell=0,1,2,\ldots$
\begin{equation}\label{eq:AffineParametersforRho=0}
\ri^+(\ell)= q \ell
\qquad
\rj^+(\ell)= p \ell
\qquad
\ri^-(\ell)= -q (\ell +1)
\qquad
\rj^-(\ell)= -p (\ell +1).
\end{equation}
Likewise, a rational critical point $\z$ with $\rho=1$ has the
following affine parameters for positive and negative chains:
for $\ell=0,1,2,\ldots$
\begin{equation}\label{eq:AffineParametersforRho=1}
\ri^+(\ell)= q (\ell +1)
\qquad
\rj^+(\ell)= p (\ell +1) -1
\qquad
\ri^-(\ell)= -q \ell
\qquad
\rj^-(\ell)= -p \ell -1.
\end{equation}
We first consider the case where $\rho=0$ and $\theta\not = 0,1$,
i.e., $\z$ is of the form $\z=(\theta,\rho)=(p/q,0)$ with $q \geqslant
2$.
In this case, the $\III$- and $\IV$-pencils are missing.
Using (\ref{eq:isigma}) and (\ref{eq:AffineParametersforRho=0}), we
see that the $\theta$-coordinate, denoted by $\theta^\I(\ell)$, of the
intersection point of the line $\rho= \ri^+(\ell) \theta- \rj^+(\ell)$, $\ell\geqslant 1$
and the positive dominant line of $\zu=(p/q,1/q)$ is given by
\begin{equation*}
\theta^\I(\ell)=
\begin{cases}
\frac{p \ell -p'}{q \ell -q'} & q' > 0\\
\frac{p \ell -p -p'}{q \ell -q -q'} & q' < 0.
\end{cases}
\end{equation*}
Similarly, the $\theta$-coordinate, denoted by $\theta^\II(\ell)$, of the
intersection point of the line $\rho= \ri^-(\ell) \theta- \rj^-(\ell)$, $\ell\geqslant 1$
and the negative dominant line of $\zu$ is given by
\begin{equation*}
\theta^\II(\ell)=
\begin{cases}
\frac{p \ell +p'}{q \ell +q'} & q' > 0\\
\frac{p \ell +p +p'}{q \ell +q +q'} & q' < 0.
\end{cases}
\end{equation*}
Using the condition for neighbouring terms in a Farey sequence, we see
that $p/q$ and $\theta^\I(\ell)$ (resp. $\theta^\II(\ell)$ and $p/q$)
are neighbours in $\mathcal{F}_{|\ri^+(\ell)|}$
(resp. $\mathcal{F}_{|\ri^-(\ell)|}$).
Thus, all the intersection points are the Farey points distinct from
$\z$ in the $\I$- and $\II$-pencils.
Since
$\mathcal{F}_{|\bi^\I|}
\subset 
\mathcal{F}_q
=
\mathcal{F}_{|\ri^+(1)|}$
and 
$\mathcal{F}_{|\bi^\II|}
\subset 
\mathcal{F}_q
\subset
\mathcal{F}_{|\ri^-(1)|}$,
all these Farey points belong to the closure of the dominant curves of $\zu$.
In the case where $\rho=1$ and $\theta\not = 0,1$, i.e., $\z$ is of
the form $\z=(\theta,\rho)=(p/q,1)$ with $q \geqslant 2$, the missing
pencils are $\I$ and $\II$.
In this case, we can do the same to show that the Farey points distinct
from $\z$ in the $\III$- and $\IV$-pencils, respectively, belong to
the positive and negative dominant curves of $\zd$.
Lastly, we consider $\z$ with $q = 1$, i.e., the four corners $(0,0)$,
$(1,0)$, $(0,1)$, $(1,1)$ of the parameter space, each of which has
only one pencil.
It is easy to see that the Farey points of the $\I$-pencil of $(0,0)$
and those of the $\II$-pencil of $(1,0)$ belong to $\rho=1$ ($0
\leqslant \theta \leqslant 1$), i.e., the negative dominant curve of
$(0,1)$ and $(1,1)$.
We can also see that the Farey points of the $\IV$-pencil of $(0,1)$ and
those of the $\III$-pencil of $(1,1)$ belong to $\rho=0$ ($0 \leqslant
\theta \leqslant 1$), i.e., the positive dominant curve of $(0,0)$ and
$(1,0)$.

We remark that the particular choice of continued fraction representation 
for $p/q$ has little effect of the above argument. It merely causes a
shift by one unity of the quantities $t^\pm$ in (\ref{eq:tpm}).
\endproof

To complete our analysis of the critical curves of the map $\rg$,
we now characterise the words of all curves incident to a 
rational point $\z$, in terms of the word at $\z$.

\begin{theorem}\label{thm:Words}
Let $\z$ be a rational critical point with denominator $q$, 
let $w^\sigma(\ell)$, $\ell\geqslant 0$ 
be the word of the $\ell$th curve adjacent to $\z$
in the quadrant $\sigma$.
Let $u^+$ (resp. $u^-$) be the word of the positive
(resp. negative) dominant curve of $\z$.
Let $v^{+}$ (resp. $v^{-}$) be the word obtained by switching
the first letter of $u^+$ (resp. $u^-$).
Then $|u^+u^-|=|u^-u^+|=q$ and
$$
w^\sigma(\ell)=\begin{cases}
u^+ \left(v^{-} u^+\right)^{\ell} &\quad \sigma=\I\\
u^- \left(u^+ v^{-}\right)^{\ell} &\quad \sigma=\II\\
u^+ \left(u^- v^{+}\right)^{\ell} &\quad \sigma=\III\\
u^- \left(v^{+} u^-\right)^{\ell} &\quad \sigma=\IV.
\end{cases}
$$
\end{theorem}
\proof
Assume first that $\rho\not=0,1$.
Let
$$
u^\sigma=\begin{cases}u^+u^-&\sigma=\I,\III\\u^-u^+&\sigma=\II,\IV.
          \end{cases}
$$
Then $u^\sigma$ is a periodic boundary word at $\z$, for the initial 
conditions $0$ or $\rho$. The periodic orbit has no other boundary 
point because $\z$ lies in the interior of the curve of both $u^+$ 
and $u^-$. Thus $|u^\sigma|=q$.
Let $w^\sigma(\ell)$ be the word of the curve adjacent to $\z$
in the quadrant $\sigma$, and let $n^\sigma(\ell)$ be the length
of this word. From (\ref{eq:ijpm}) we find
$$
n^\sigma(\ell)=|\ri^\pm(\ell)|=q(\{\pm\tau\}+\ell),
$$
with the usual convention on sign, and the quotient of division
of $n^\sigma(\ell)$ by $q$ is given by
$$
\lfloor n^\sigma(\ell)/q\rfloor=\lfloor \{\pm\tau\}+\ell\rfloor =\ell.
$$
The word $w^\sigma(\ell)$ is now computed from theorem 
\ref{thm:CriticalCurves}, part ii). 
Then $w^\sigma(\ell)$ will consist of $\ell$ repetitions of a 
modification of $u^\sigma$,
followed by a modification of $u^\pm$,
where the modifications are performed on the symbols 
congruent to $|\ri^\sigma(0)|$ modulo $q$ if the curve is on the right of $\z$
($\sigma=\I,\IV$), or those congruent to $0$
modulo $q$, except the first symbol,
if the curve is on the left of $\z$ ($\sigma=\II,\III$).
This gives $w^\sigma(\ell)$ in the statement.

We now consider the case $\rho=0$.
By definition, $u^+ = \varepsilon$ and $u^- = b^q$, where $\varepsilon$
denotes the empty word [see discussion preceding (\ref{eq:zud})].
Obviously, $|u^+u^-|=|u^-u^+|=q$.
For $\sigma=\I,\II$, we denote by $w^\sigma(\ell)$ the word of the curve
adjacent to $\z$ in the quadrant $\sigma$ and put $n^\sigma(\ell) =
|w^\sigma(\ell)|$, as above.
From (\ref{eq:AffineParametersforRho=0}) we have
$$
n^\I(\ell)=|\ri^+(\ell)|= q \ell
\hskip 30pt
n^\II(\ell)=|\ri^-(\ell)|= q (\ell +1).
$$
The word $w$ of length $n^\I(\ell)=q \ell$ at $\z$ is given by $w=b^{q
  \ell}$.
Thus, theorem \ref{thm:CriticalCurves} ii) gives
$$
w^\I(\ell) = \left(a b^{q-1}\right)^\ell = \varepsilon \left(a b^{q-1}\varepsilon\right)^\ell
= u^+ \left(v^{-} u^+\right)^{\ell}.
$$
Likewise, the word $w$ of length $n^\II(\ell)= q (\ell +1)$ at $\z$ is
$w=b^{q (\ell +1)}$, and the same theorem gives
$$
w^\II(\ell) = b^{q} \left(a b^{q-1}\right)^\ell = b^{q} \left(\varepsilon a b^{q-1}\right)^\ell = u^- \left(u^+ v^{-}\right)^{\ell}.
$$
This completes the proof of the case $\rho=0$.
We can do the same for the proof of the case $\rho=1$.
\endproof

\section{Triple points and convergents}\label{section:TriplePoints}

\begin{figure}[t]
\hspace*{-70pt}
\begin{minipage}{5cm}
\centering
      \includegraphics[scale=0.60,clip]{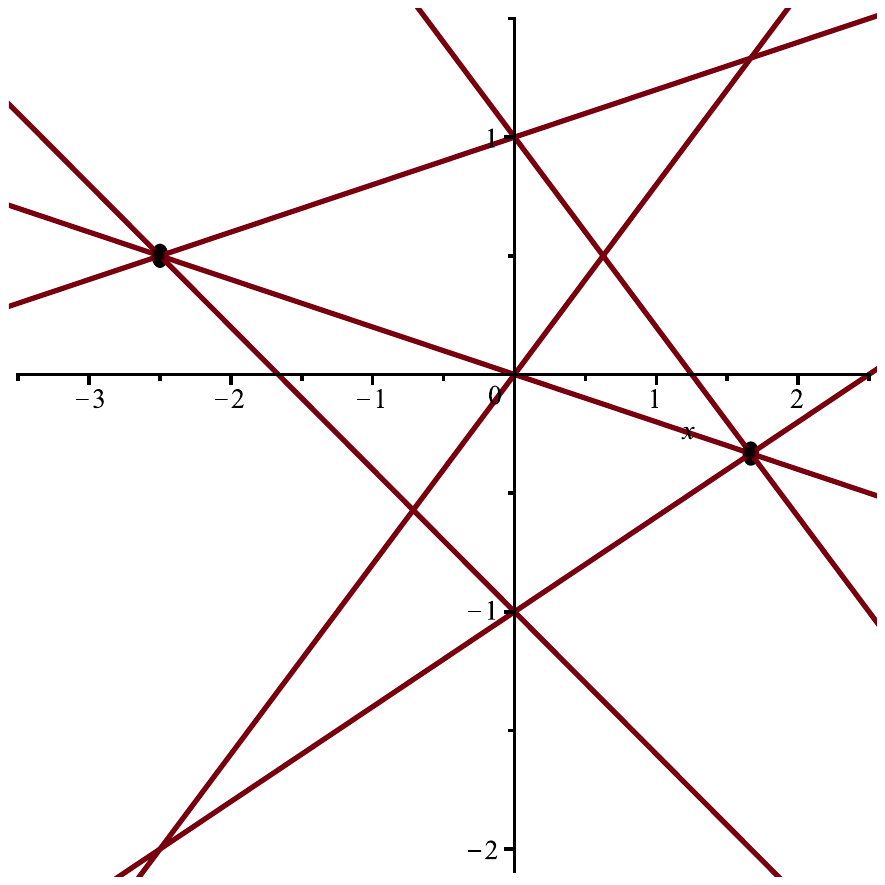}
\end{minipage}
\hspace*{40pt}
\begin{minipage}{5cm}
\centering
      \includegraphics[scale=0.60,clip]{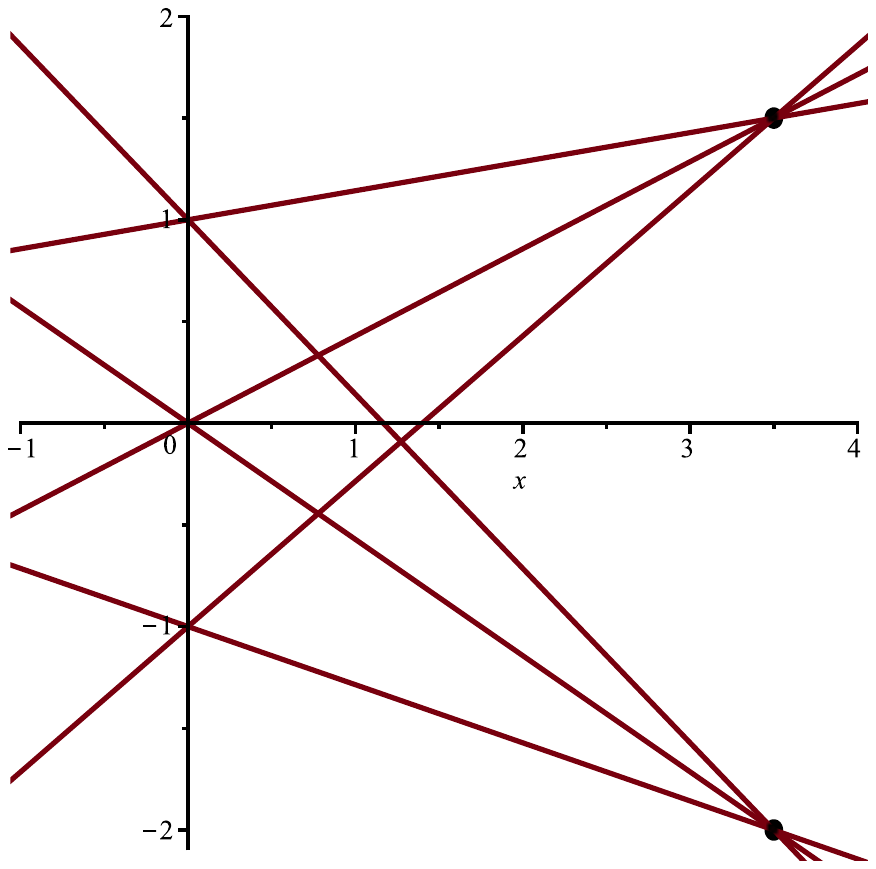}
\end{minipage}
\vspace*{-270pt}
\caption{\label{fig:TriplePoints} 
Examples of triple points for $\z=(3/5,2/5)$ (left, type I)
and for $\z=(3/7,2/7)$ (right, type II), as specified in 
theorem \ref{thm:TriplePoints}. 
The coordinates are normalised by shifting $\z$ to the origin, 
and then scaling $\rho$ by $q$, and $\theta$ by $q^2$.
}
\end{figure}

In the previous section we considered the solutions of 
(\ref{eq:CriticalPoint}) for fixed rational $(\theta,\rho)$, resulting 
in infinitely many chains passing through that point.
A global view on symbolic dynamics may be gained by considering
a different set of chains, namely those corresponding to
the set $\cN_n$ of solutions of (\ref{eq:CriticalPoint}) 
with $|\ri|$ bounded by $n$, and $\rj$ subject to the bounds 
$J(\ri)$ given in (\ref{eq:Jbounds}):
\begin{equation}\label{eq:Net}
\cN_n=\bigcup_{|\ri|\leqslant n\atop J(\ri)} \cLij.
\end{equation}
By construction, the boundary words associated to the chains in
$\cN_n$ contain all possible critical words of length not
exceeding $n$ and either sign.
For brevity, we shall not develop the analysis of $\cN_n$ here,
but merely prove a geometric theorem, illustrated in figure 
\ref{fig:TriplePoints}, which deals with a subset of $\cN_n$ consisting
of six chains near a rational critical point. 
This result displays a connection between intersections of 
chains and convergents of continued fractions.

Let $\z=(p/q,r/s)$ be a rational critical point with two neighbours
$\zu$ and $\zd$, that is, $r/s\not=0,1$.
We recall (see beginning of section \ref{section:RationalPoints})
that $p_k/q_k$ denote the convergents of the continued fraction 
expansion of $p/q=p_n/q_n$, chosen so to have the last coefficient 
equal to one.
Then $p_{n-1}/q_{n-1}$ is the rational closest to $p/q$ among 
those with denominator less than $q$. (This is not the case if the
last coefficient is greater than one.)
We define
\begin{equation}\label{eq:mupsi}
\mu=2\flo{\tau}-\flo{\tau^-}-\flo{\tau^+}
\hskip 30pt
\psi_\mu(x)=\begin{cases} \lfloor x\rfloor & \mu=+1\\
                          \lceil  x\rceil  & \mu=-1,
            \end{cases}
\end{equation}
where $\tau,\tau^-,\tau^+$ are given in (\ref{eq:tau}).

A \textbf{triple point} of $\z$ is a common point of 
three concurrent dominant lines, which comprise one line from 
each of $\z,\zu,\zd$. The next result characterises all
triple points.

\begin{theorem}\label{thm:TriplePoints}
Let $\z=(p/q,r/s)$ be a rational critical point with 
two neighbours, let $\mu$ and $\psi$ be given by (\ref{eq:mupsi}), and let
$$
\chi^{(1)}_\mu=\biggl(\frac{p_{n-1}}{q_{n-1}},
  \frac{\psi_\mu(\tau)}{q_{n-1}}(-1)^{n-1}\biggr)
\hskip 35pt
\chi^{(2)}_\mu=\biggl(\frac{p_{n-2}}{q_{n-2}},
    \frac{rq_n}{sq_{n-2}}
    -\frac{\psi_\mu(\tau)}{q_{n-2}}(-1)^{n-1}\biggr).
$$
Then all $\z$ have two triple points, specified below.\\
The triple points of $\z=(p/q,1/q)$ are:
\vspace*{-10pt}
\begin{enumerate}
\item [] $\chi^{(1)}_{-1}$ and $\chi^{(2)}_{-1}$ \/ if $n$ is odd \/ 
            (type I);
\item [] $\chi^{(2)}_{+1}$ and $\chi^{(2)}_{-1}$ \/ if $n$ is even \/ 
            (type II).
\end{enumerate}
The triple points of $\z=(p/q,(q-1)/q)$ are:
\vspace*{-10pt}
\begin{enumerate}
\item [] $\chi^{(1)}_{+1}$ and $\chi^{(2)}_{+1}$ \/ if $n$ is odd \/ 
            (type I);
\item [] $\chi^{(2)}_{+1}$ and $\chi^{(2)}_{-1}$ \/ if $n$ is even \/ 
            (type II).
\end{enumerate}
The triple points of all other $\z$ are:
\vspace*{-10pt}
\begin{enumerate}
\item [] $\chi^{(1)}_\mu$ and $\chi^{(2)}_\mu$ \/ if $\mu\not=0$ \/ 
            (type I);
\item [] $\chi^{(2)}_{+1}$ and $\chi^{(2)}_{-1}$ \/ if $\mu=0$ \/ 
            (type II).
\end{enumerate}
\end{theorem}
The cases listed above as I and II will be referred to as the 
\textit{generic cases} ---see figure \ref{fig:TriplePoints}.

\proof
Assume first that $r/s \not= 1/q, (q-1)/q$.
\def\sd{{\mu_1}}%
\def\s0{{\mu_2}}%
\def\su{{\mu_3}}%
We choose a dominant line from each point, and form the matrix whose
rows represent the chosen lines.
We will show that the determinant of exactly two such matrices
must vanish. These determinants take the form [cf.~eq.(\ref{eq:ijpm})]
$$
D(\sd,\s0,\su)=
\left\vert
\begin{matrix}
\sd\{\sd\tau^-\} & -1 & 1\\
\s0\{\s0\tau\} & 0 & 1\\
\su\{\su\tau^+\} &1 & 1
\end{matrix}
\right\vert
\qquad
\sd,\s0,\su\in\{+1,-1\},
$$
where the rows represent the lines through $\zd,\z$ and $\zu$,
respectively, and for each row there is one choice of sign, 
to select the positive or negative line through that point.
To normalise the determinant, we have placed the origin at $\z$, 
divided the first column by $q$, and multiplied the second column by $q$.

Using the identities $\tau^+-2\tau+\tau^-=0$,
$\{x\}=x-\lfloor x\rfloor$, and $-\lfloor -x\rfloor=\lceil x\rceil$,
we obtain
\begin{eqnarray*}
D(\sd,\s0,\su)
&=& -\sd\{\sd\tau^-\}+2\s0\{\s0\tau\}-\su\{\su\tau^+\}\\
&=& \sd\lfloor\sd\tau^-\rfloor-2\s0\lfloor\s0\tau\rfloor
  +\su\lfloor\su\tau^+\rfloor\\
&=& \psi_\sd(\tau^-)-2\psi_\s0(\tau)+\psi_\su(\tau^+),
\end{eqnarray*}
where $\psi$ is defined in (\ref{eq:mupsi}).
Thus all eight determinants are integer, and the
vanishing of the determinant is equivalent to $\psi_{\mu_2}(\tau)$
being the mid-point of $\psi_{\sd}(\tau^-)$ and
$\psi_{\su}(\tau^+)$, that is, that 
$\psi_\sd(\tau^-),\psi_\s0(\tau),\psi_\su(\tau^+)$ form 
an arithmetic progression.

The sequence 
$\lfloor \tau^-\rfloor$, $\lfloor \tau\rfloor$, $\lfloor \tau^+\rfloor$
is non-decreasing if $q'>0$ and non-increasing if $q'<0$.
Since $|q'|/q<1$, we have 
$|\lfloor \tau^\pm\rfloor-\lfloor\tau\rfloor|\leqslant 1$,
and hence such a sequence can assume only one of the following 
values
\begin{equation}\label{eq:floorseq}
(k,k,k),\quad(k,k,k\pm1),\quad(k,k\pm1,k\pm1),\quad(k,k\pm1,k\pm2)
\end{equation}
for some integer $k$, where the sign is that of $q'$.

Since $r/s \not= 1/q, (q-1)/q$, 
none of the $\tau$s is an integer, and so 
we have 
$$
\psi_{-1}(x)=\lceil x\rceil=\lfloor x\rfloor+1=\psi_{1}(x)+1,
      \hskip 30pt x\in\{\tau^-,\tau,\tau^+\}.
$$
With this in mind, we match each sequence in (\ref{eq:floorseq})
with two $\mu$-sequences so as to transform it into an arithmetic
progression.
\begin{equation}\label{eq:mu-words}
\begin{array}{ccccccccc}
 &     &\lfloor\tau^-\rfloor&\lfloor\tau\rfloor&\lfloor\tau^+\rfloor&
       &\sd\s0\su&&\sd\s0\su\\
\noalign{\vskip 15pt}
1&\quad&k&k&k&\qquad&+++&\quad&---\\
2&\quad&k&k+1&k+2&\qquad&+++&\quad&---\\
3&\quad&k&k-1&k-2&\qquad&+++&\quad&---\\
\noalign{\vskip 3pt}
4&\quad&k&k&k+1&\qquad&--+&\quad&+--\\
5&\quad&k&k+1&k+1&\qquad&-++&\quad&++-\\
6&\quad&k&k&k-1&\qquad&++-&\quad&-++\\
7&\quad&k&k-1&k-1&\qquad&+--&\quad&--+
\end{array}
\end{equation}

The above argument shows that in general a rational critical point
has two $\mu$-sequences corresponding to two triple points, which
we now compute. 
We recall that $p_k/q_k$ are the convergents of $p/q=p_n/q_n$.
We let $p'=(-1)^{n-1}p_{n-1}$ and $q'=(-1)^{n-1}q_{n-1}$ 
[cf.~(\ref{eq:tau})], whence $pq'-p'q=1$.

We choose to consider the intersection of the line through $\zd$ with
slope $q\sd\{\sd\tau^-\}$ with the line through $\z$ with slope
$q\s0\{\s0\tau\}$. 
Letting
\begin{equation}\label{eq:NablaDefinition}
\nabla=\sd\lfloor\sd\tau^-\rfloor-\s0\lfloor\s0\tau\rfloor
 =\psi_{\mu_1}(\tau^-)-\psi_{\mu_2}(\tau),
\end{equation}
we obtain for the rotation number
\begin{eqnarray*}
\theta&=&\frac{p}{q}+\frac{1}{q^2(\sd\{\sd\tau^-\}-\s0\{\s0\tau\})}
 =\frac{pq'+pq\nabla-1}{q(q'+q\nabla)}=\frac{p'+p\nabla}{q'+q\nabla}\\
 &=&\frac{(-1)^{n-1}p_{n-1}+p_n\nabla}{(-1)^{n-1}q_{n-1}+q_n\nabla}.
\end{eqnarray*}
Having assumed that the last coefficient of the continued fractions
of $p_n/q_n$ is unity, we have $p_n-p_{n-1}=p_{n-2}$ and
$q_n-q_{n-1}=q_{n-2}$. 
We find
\begin{equation}\label{eq:ThetaTriplePoints}
\theta=\begin{cases} {p_{n-1}}/{q_{n-1}} & \quad\nabla=0\\
                     {p_{n-2}}/{q_{n-2}} & \quad\nabla=(-1)^n.\\
       \end{cases}
\end{equation}
We now show that $\nabla$ does not assume any other value.
With reference to table (\ref{eq:mu-words}) we find:
\begin{equation}\label{eq:nabla}
\begin{array}{cccclclcl}
\sd\s0&\quad&\nabla&\quad&\nabla=0&\quad&\nabla=-1&\quad&\nabla=+1\\
\noalign{\vskip 15pt}
++&&\flo{\tau^-}-\flo{\tau}&&1,6&&2,5&&3\\
\noalign{\vskip 3pt}
--&&\cei{\tau^-}-\cei{\tau}&&1,4&&2  &&3,7\\
\noalign{\vskip 3pt}
+-&&\flo{\tau^-}-\cei{\tau}&&7  &&4  &&\\
\noalign{\vskip 3pt}
-+&&\cei{\tau^-}-\flo{\tau}&&5  &&   &&6\\
\end{array}
\end{equation}
where the last three columns list the cases in (\ref{eq:mu-words})
where each value of $\nabla$ is attained.
Considering that $n$ is odd in cases 2,4,5 and $n$ is even in cases
3,6,7, we see that, if $\nabla\not=0$, then
$\nabla$ is given by $(-1)^n$.

We now divide the rows of table (\ref{eq:mu-words}) into two groups,
namely the rows 1--3 where $\mu=0$ [$\mu$ is defined in (\ref{eq:mupsi})]
(that is, $\flo{\tau^-},\flo{\tau},\flo{\tau^+}$ form
an arithmetic progression), and the rows 4--7 where $\mu\not=0$.
In the latter cases we verify that both values of $\theta$ listed
in (\ref{eq:ThetaTriplePoints}) occur. Thus the triple points are 
located on opposite sides of $\z$. 

In the former case only one value of $\theta$ occur, and both triple 
points lie on the same side of $\z$. To determine the value of 
$\theta$, we note  that case 1 in (\ref{eq:mu-words}) does not 
actually occur.
Indeed this would require 
$$
|\tau^+-\tau^-|=2\frac{q_{n-1}}{q_n} < 1
$$
or $q_{n-1}/q_n < 1/2$. 
Since $q_{n-2} + q_{n-1}= q_n$, this would give $q_{n-2} > q_{n-1}$,
which is a contradiction.
Therefore $\theta=p_{n-2}/q_{n-2}$ for both points.

It remains to determine the value of $\rho$ at the triple point.
The dominant lines of $\z$ have equation
$$
\rho-\frac{r}{s}=q_n\mu_2\{\mu_2\tau\}\Bigl(\theta-\frac{p_n}{q_n}\Bigr),
$$
where $\mu_2$ is determined according to table (\ref{eq:mu-words}).
Letting $\theta=p_{n-1}/q_{n-1}$, we find
\begin{eqnarray*}
\rho&=&\frac{r}{s}+q_n\mu_2(\mu_2\tau-\flo{\mu_2\tau})
    \frac{(-1)^n}{q_{n-1}q_n}\\
&=&-\mu_2\flo{\mu_2\tau}\frac{(-1)^n}{q_{n-1}}
   =\psi_{\mu_2}(\tau)\frac{(-1)^{n-1}}{q_{n-1}},
\end{eqnarray*}
and we have obtained the triple point $\chi^{(1)}_\mu$.
One verifies that $\mu_2=\mu$.
Letting $\theta=p_{n-2}/q_{n-2}$, we find, using \cite[theorem 151]{HardyWright}
$$
\rho=\frac{rq_n}{sq_{n-2}}-\psi_{\mu_2}(\tau)\frac{(-1)^{n-1}}{q_{n-2}},
$$
which is $\chi^{(2)}_\mu$. In this case, from table (\ref{eq:mu-words}),
we see that if $\mu\not=0$, then $\mu_2=\mu$, while for $\mu=0$,
both signs are allowed.

We now discuss the cases $r/s=1/q$ and $r/s=(q-1)/q$.
If $q=2$, then $\z=(1/2,1/2)$, and we see immediately that the triple
points of $\z$ are $(0,0)$ and $(0,1)$ [see discussion preceding
  (\ref{eq:zud})].

Let $r/s=(q-1)/q$ with $q \geqslant 3$.
Then, $\zu$ is on $\rho=1$, and the slope of a dominant line of $\zu$
is $q (\mu_3 +1)/2$ from (\ref{eq:AffineParametersforRho=1}).
This gives the condition for three concurrent dominant lines as
follows:
\begin{equation}\label{eq:ThreeConcurrentDominantLinesFor(q-1)/q}
\psi_\sd(\tau^-)-2\psi_\s0(\tau)+\tau^+ - \frac{\mu_3 +1}{2} =0.
\end{equation}
As noted above, case 1 in (\ref{eq:mu-words}) does not occur.
Other than case 1, cases 3, 4, 5, 6 do not occur, too:
In fact, cases 3, 5, 6 do not occur since $\tau^+ = q'$ is an integer
and $0 < |\tau^+ - \tau| = q_{n-1}/q < 1$ for $q \geqslant 2$.
Likewise, case 4 does not occur since $|\tau^+ - \tau^-| = 2 q_{n-1}/q
> 1$ for $q \geqslant 3$.
For remaining cases 2, 7, we give $\mu$-sequences satisfying
(\ref{eq:ThreeConcurrentDominantLinesFor(q-1)/q}):
$$
\begin{array}{ccccccccc}
 &     &\lfloor\tau^-\rfloor&\lfloor\tau\rfloor&\lfloor\tau^+\rfloor = \tau^+&
       &\sd\s0\su&&\sd\s0\su\\
\noalign{\vskip 15pt}
2&\quad&k&k+1&k+2&\qquad&-++&\quad&++-\\
7&\quad&k&k-1&k-1&\qquad&+++&\quad&---
\end{array}
$$
In case 2, $n$ is odd, and $\nabla$ defined in
(\ref{eq:NablaDefinition}) is $0$ if $\sd\s0\su = -++$, and $-1 =
(-1)^n$ if $\sd\s0\su = ++-$.
Thus, (\ref{eq:ThetaTriplePoints}) also specifies the $\theta$ value
of a common point for this case.
Since $\s0 = +1$ for the both $\mu$-sequences, we have
$\chi^{(1)}_{+1}$ and $\chi^{(2)}_{+1}$ as the triple points for case
2.
In case 7, $n$ is even, and $\nabla = 1 = (-1)^n$ for the both
$\mu$-sequences.
Thus, by (\ref{eq:ThetaTriplePoints}), only
$\theta={p_{n-2}}/{q_{n-2}}$ occurs, and we have $\chi^{(2)}_{+1}$ and
$\chi^{(2)}_{-1}$ as the triple points for case 7.

In the case of $r/s=1/q$ with $q \geqslant 3$, $\zd$ is on $\rho=0$,
and the slope of a dominant line of $\zd$ is $q (\mu_1 -1)/2$
[cf.~(\ref{eq:AffineParametersforRho=0})].
This gives the following condition for three concurrent dominant lines:
\begin{equation}\label{eq:ThreeConcurrentDominantLinesFor1/q}
\tau^- - \frac{\mu_1 -1}{2} -2\psi_\s0(\tau)+\psi_\su(\tau^+) =0.
\end{equation}
Incidentally, $\tau^- = 0$.
By an argument similar to that given above, we can see that only cases
3, 4 in (\ref{eq:mu-words}) occur.
By (\ref{eq:ThreeConcurrentDominantLinesFor1/q}), we have the
following table:
$$
\begin{array}{ccccccccc}
 &     &\lfloor\tau^-\rfloor = \tau^- &\lfloor\tau\rfloor&\lfloor\tau^+\rfloor&
       &\sd\s0\su&&\sd\s0\su\\
\noalign{\vskip 15pt}
3&\quad&k&k-1&k-2&\qquad&+++&\quad&---\\
4&\quad&k&k&k+1&\qquad&--+&\quad&+--
\end{array}
$$
where $k$ is actually 0.
Let
$$
\nabla'=\su\lfloor\su\tau^+\rfloor-\s0\lfloor\s0\tau\rfloor
 =\psi_{\mu_3}(\tau^+)-\psi_{\mu_2}(\tau).
$$
Considering the intersection of dominant lines through $\zu$ and $\z$,
we have the $\theta$ value of a common point as follows:
$$
\theta=\begin{cases} {p_{n-1}}/{q_{n-1}} & \quad\nabla'=0\\
                     {p_{n-2}}/{q_{n-2}} & \quad\nabla'=(-1)^{n-1}.\\
       \end{cases}
$$
In case 3, $n$ is even, and $\nabla' = -1 = (-1)^{n-1}$ for the both
$\mu$-sequences.
Thus, we have $\chi^{(2)}_{+1}$ and $\chi^{(2)}_{-1}$ as the triple
points for case 3.
In case 4, $n$ is odd, and $\nabla'=0$ if $\sd\s0\su = --+$, and
$\nabla'= 1 = (-1)^{n-1}$ if $\sd\s0\su = +--$.
Since $\s0 = -1$ for the both $\mu$-sequences, we have
$\chi^{(1)}_{-1}$ and $\chi^{(2)}_{-1}$ as the triple points for case
4.

We note that $\z=(1/2,1/2)$ can be classified into either
$\z=(p/q,1/q)$ or $(p/q,(q-1)/q)$ since $n = 2$ is even and its triple
points $(0,0)$ and $(0,1)$ are $\chi^{(2)}_{+1}$ and
$\chi^{(2)}_{-1}$.

The above argument applies to continued fraction representations
for which the last coefficient is chosen to be unity. 
If instead the last coefficient is greater than 1, then 
we write $p/q=\bar p_{\bar n}/\bar q_{\bar n}$, where $\bar n=n-1$.
We find
$$
\bar q_{\bar n-1}=q_{n-2}
\qquad\mbox{and}\qquad
\bar q_{\bar n}-\bar q_{\bar n-1}=q_{n-1}
$$
and similarly for $\bar p_{\bar n}$. 
It follows that the two values of $\theta$ in 
(\ref{eq:ThetaTriplePoints}) are merely exchanged, 
yielding the same result.
\endproof

Considering the slopes of the lines at a triple point,
from the above theorem and lemma \ref{lma:Endpoints} iii) 
we obtain at once:

\begin{corollary} \label{cor:TriplePoints}
A triple point of type I is the Farey point of at least one of the three
concurrent curves; for type II, it is the Farey point of at least two curves.
\end{corollary}

We finally consider sequences of rational critical points approaching 
an irrational one. 
Recall (see beginning of section \ref{section:BasicProperties})
that the point $\z=(\theta,\rho)$ is a critical point iff $\rho$
belongs to the doubly infinite orbit through the origin.
Therefore $\z$ has the form 
\begin{equation}\label{eq:zeta}
\z=(\theta,\{\ri\theta\})=(\theta,\ri\theta-\rj)\qquad 
 \rj=\lfloor \ri\theta\rfloor,
\end{equation}
for some $\ri\in\Z$. However, the continued fractions of $\theta$ 
and $\ri\theta$ will in general be unrelated, and so it is not 
possible to construct a sequence of rational critical points 
that consists of convergents of both $\theta$ and $\rho$.
We shall therefore prioritise the convergents of $\theta$,
and select rational approximants for $\rho$ so that the affine
parameters $(\ri,\rj)$ remain the 
same throughout the approximation.

Thus choose $\theta\not\in\Q$ and $\ri$ in (\ref{eq:zeta}), and 
consider the sequence of convergents $p_k/q_k, k\geqslant 0$ of
the continued fraction expansion of $\theta$.
From (\ref{eq:zeta}) we obtain the following sequence of 
rational critical points
$$
\z_k=\Bigl(\frac{p_k}{q_k},\frac{\ri p_k-\rj q_k}{q_k}\Bigr)\qquad k\geqslant 0,
$$
with $\z_k\to\z$ as $k\to\infty$.

By construction, all points $\z_k$ are collinear.
Now, the dominant affine parameters $(\bi^\pm,\bj^\pm)$ at $\z_k$
have the property that $|\bi^+-\bi^-|=q_k$ (see beginning of section
\ref{section:RationalPoints}).
Since the pair $(\ri,\rj)$ is fixed and $q_k\to\infty$, it follows that 
$(\ri,\rj)$ is a pair of dominant affine parameters of the point $\z_k$ 
for all sufficiently large $k$, while the components of the other 
dominant pair diverge to infinity. 

\end{document}

%% file: fig-Chain3.tex
\begin{tikzpicture}[scale=70]
\tikzmath{\xm=0.02;}
\tikzmath{\xs=0.01;}
\tikzmath{\xxs=0.005;}

\draw[black, thin] (5/7,0) -- (6/7,0);

\filldraw[blue] (5/7,0) circle (0.06pt);
\node at (5/7,-\xs) {$\theta^-=\frac{5}{7}$};
\node at (5/7,-\xs-\xs) {$b^7$};
\node at (5/7,-\xs-\xs-\xs) {$\varepsilon$};

\node at (5/7+\xm,\xs) {$ab^6$};

\filldraw[blue] (3/4,0) circle (0.06pt);
\node at (3/4,-\xs) {$\frac{3}{4}$};
\node at (3/4,-\xs-\xs) {$ab^3ab^2$};
\node at (3/4,-\xs-\xs-\xs) {$ab^2$};

\node at (3/4+\xm+\xxs,\xs) {$ab^2a^2b^2$};

\filldraw[blue] (4/5,0) circle (0.06pt);
\node at (4/5,-\xs) {$\frac{4}{5}$};
\node at (4/5,-\xs-\xs) {$ab^2a^3b$};
\node at (4/5,-\xs-\xs-\xs) {$ab$};

\node at (4/5+\xs+\xxs,\xs) {$aba^4b$};

\filldraw[blue] (5/6,0) circle (0.06pt);
\node at (5/6,-\xs) {$\frac{5}{6}$};
\node at (5/6,-\xs-\xs) {$aba^5$};
\node at (5/6,-\xs-\xs-\xs) {$a$};

\node at (5/6+\xs+\xxs,\xs) {$a^7$};

\filldraw[blue] (6/7,0) circle (0.06pt);
\draw (6/7,0) circle (0.06pt);
\node at (6/7,-\xs) {$\theta^+=\frac{6}{7}$};
\node at (6/7,-\xs-\xs) {$a^7$};
\node at (6/7,-\xs-\xs-\xs) {$\varepsilon$};

\end{tikzpicture}